\documentclass{nyj}

\title{RIEMANN HYPOTHESIS AND SUPERCONFORMAL INVARIANCE} 

\author{Matti Pitk\"anen}  
\address{Dept. of Physics, University of Helsinki, Helsinki, Finland.} 
\email{matpitka@rock.helsinki.fi, URL: http://www.physics.helsinki.fi/$\tilde{~}$matpitka/}   

\thanks{
I am grateful for my son Timo for stimulating discussions relating
to  p-adic numbers  in the context
of the Riemann Zeta function.  
I would like  to thank Carlos Castro for inspiring
email conversations relating to Riemann Zeta and for turning my attention to Hilbert-Polya hypothesis and the possible physics involved with it. 
I also want  to   express
my deep gratitude to  Matthew Watkins for providing me
with information about Riemann Zeta and for generous help, in particular for reading the manuscript and pointing out several inaccuracies and errors.}


\subjclass{Number Theory}

\newcommand{\vm}{\vspace{0.2cm}}
\newcommand{\vl}{\vspace{0.4cm}}
\newcommand{\per}{\hspace{.2cm}}

\begin{document} 

\begin{abstract}  
 A strategy for proving (not a proof of, as was the first over-optimistic belief) 
the Riemann hypothesis is suggested. The vanishing of Riemann Zeta reduces to an orthogonality condition for the eigenfunctions of a non-Hermitian operator $D^+$ having the zeros of  Riemann Zeta as its eigenvalues. The construction of $D^+$ is inspired by the  conviction  that Riemann Zeta is associated with a physical system allowing superconformal transformations as  its  symmetries and second quantization in terms of the representations of superconformal algebra. The eigenfunctions  
of $D^+$ are analogous to the so called  coherent states and in general not orthogonal  to each other. The states orthogonal to a vacuum state (having a negative norm squared) correspond  to the zeros of Riemann Zeta. The 
physical states  having a positive norm squared correspond to the zeros of Riemann Zeta at the critical line. Riemann hypothesis follows by  reductio ad absurdum from the hypothesis that  ordinary superconformal algebra acts as gauge  symmetries for all  coherent states orthogonal to  the  vacuum state, including  also the non-physical might-be coherent states  off from the critical line.
 \end{abstract} 
\maketitle



\section{Introduction}

The Riemann hypothesis \cite{Riemann, Titchmarch}
 states that the non-trivial  zeros (as opposed
to zeros at $s= -2n$, $n\geq 1$ integer) of Riemann Zeta function obtained by
analytically continuing the function 

\begin{eqnarray}
\zeta (s) = \sum_{n=1}^{\infty} \frac{1}{n^s} 
\end{eqnarray}

\noindent from the region $Re[s]>1$
to the  entire complex plane,  lie on the line  $Re[s]=1/2$. 
Hilbert and Polya \cite{Hpolya} conjectured a long time ago
that the non-trivial zeroes of Riemann Zeta function could have  spectral
interpretation in terms of the eigenvalues of a suitable self-adjoint differential operator $H$ such that the eigenvalues of
this operator correspond to the imaginary
parts  of   the  nontrivial zeros $z=x+iy$
of $\zeta$. One can however consider a variant of this  hypothesis
stating that  the eigenvalue
spectrum of a   non-Hermitian operator   $D^+$
contains  the non-trivial zeros of $\zeta$. 
The eigenstates in question are 
eigenstates of an annihilation operator type operator $D^+$
and analogous to the  so called coherent states encountered in 
quantum physics \cite{field}. In particular, the eigenfunctions 
are in general non-orthogonal and this is a quintessential element of the
the proposed strategy of proof.

In the following an explicit operator having as its
eigenvalues the non-trivial zeros of $\zeta$ is constructed. 

a) The construction relies
crucially on the interpretation of the vanishing of
$\zeta$  as an orthogonality condition in a  Hermitian metric which
is is a priori more general than Hilbert space inner product.

b) Second basic element is  the scaling invariance motivated by the belief
that $\zeta$ is associated with a physical
system which has superconformal transformations \cite{sconf} as its symmetries.
This vision  was inspired by  the  generalization of $\zeta$ 
and the Riemann hypothesis  to a $p$-adic context forcing the 
sharpening of the Riemann hypothesis
to the conjecture that $p^{iy}$ defines a
rational phase factor for all non-trivial zeros $x+iy$ of $\zeta$ and 
for all primes $p$
\cite{octo}. Here however only the Riemann hypothesis is discussed. 

The core elements of the construction are following.

a)  All complex
numbers are  candidates for the eigenvalues of $D^+$
and   genuine
eigenvalues are selected by the requirement that the
condition  $D^{\dagger}=D^+$
holds true in the set of the genuine eigenfunctions. This condition
is equivalent with the Hermiticity of the Hermitian metric defined
by a function proportional to $\zeta$.

b)  The  eigenvalues turn out to consist of $z=0$ and the
non-trivial zeros of $\zeta$ and only the zeros  with $Re[z]=1/2$
correspond to  the eigenfunctions  having  real  norm.
The vanishing of $\zeta$  tells that the 'physical'  positive norm
eigenfunctions, which are {\it not} orthogonal to each other,   are
orthogonal to the the  'unphysical'  negative norm  eigenfunction
associated with the eigenvalue $z=0$.   The requirement that
the Hermitian form in question defines an inner product implies
that the the sums $z_{12}=1+i(y_1+y_2)$  of 
the zeros $z=1/2+y_i$, $i=1,2$,   of $\zeta$  correspond to almost-zeros
of $\zeta$ for large values of $y_1+y_2$.

c)  The theory allows supersymmetrization and second quantization
in tems of the representations of a superconformal algebra associated
with the operator $D^+$ and containing the ordinary superconformal algebra 
\cite{sconf} as its subalgebra.  The states on the critical line
correspond to the  representations of the ordinary superconformal
algebra acting as gauge symmetries.  If one requires that
this is also the case for  the might-exist unphysical coherent states 
orthogonal  to the vacuum state but off from the critical line, Riemann 
hypothesis  follows by a reductio ad absurdum argument.

\section{Modified form of the Hilbert-Polya conjecture}

 One can  modify  the Hilbert-Polya  conjecture by assuming 
scaling invariance and giving up
the Hermiticity of the Hilber-Polya operator. This means
introduction of the    non-Hermitian operators $D^+$ and $D$
which are Hermitian conjugates of each other such that
 $D^{+}$ has the nontrivial zeros of $\zeta$
as its  complex eigenvalues

\begin{eqnarray}
\begin{array}{ll}
D^+ \Psi= z \Psi \per .\\ 
\end{array}
\end{eqnarray}

\noindent The counterparts of the  so called  coherent states
\cite{field} are in question and 
the eigenfunctions  of $D^+$ are not expected to be orthogonal in general.
  The following construction
is based on the idea that
$D^{+}$  also allows  the eigenvalue $z=0$ and 
that the vanishing of $\zeta$ at $z$ expresses the  orthogonality
of the states with eigenvalue $z=x+iy\neq 0$ and the state
with eigenvalue  $z=0$ which turns out to have a negative norm.

  The trial 

\begin{eqnarray}
\begin{array}{ll}
D  = L_0+ V  \per ,  &D^+ = -L_0+V\\
\\
L_0=t\frac{d}{dt}\per ,& V= \frac{dlog(F)}{d(log(t))} =t\frac{dF}{dt}\frac{1}{F} \per \\ 
\end{array}
\end{eqnarray}

\noindent is motivated by the requirement of   invariance
with respect to scalings $t\rightarrow \lambda t$ and
$F\rightarrow \lambda F$. The range of variation for
the variable $t$ consists of non-negative real numbers $t\geq 0$. 
The scaling invariance implying
conformal invariance  (Virasoro generator $L_0$ represents
scaling which plays a  fundamental role in the
superconformal theories \cite{sconf})
is motivated by the  belief that $\zeta$ codes for  the physics
of a quantum critical system having, not only supersymmetries \cite{Berry}, 
but  also superconformal transformations as its basic symmetries
 \cite{octo,Castro}.

\section{Formal solution of the eigenvalue equation for operator $D^{+}$}

  One can formally solve the eigenvalue equation 

\begin{eqnarray}
D^+ \Psi_z= \left[-t\frac{d}{dt} +t\frac{dF}{dt}\frac{1}{F}\right]\Psi_z = 
 z \Psi_z \per .
\end{eqnarray}

\noindent for  $D^{+}$ by  factoring the eigenfunction to a product: 

\begin{eqnarray}
\Psi_z= f_z F\per .
\end{eqnarray}

\noindent The substitution into the eigenvalue equation gives

\begin{eqnarray}
L_0 f_z= t\frac{d}{dt}f_z = -zf_z
\end{eqnarray}

\noindent allowing as its solution the functions

\begin{eqnarray}
f_z(t) = t^{z} \per .
\end{eqnarray}

\noindent These functions are nothing but eigenfunctions of the scaling
operator $L_0$ of the superconformal
algebra analogous to the eigenstates of a translation operator.
 A priori all complex numbers
$z$ are candidates for the eigenvalues of $D^+$  and one must select
the genuine eigenvalues by applying  the requirement
$D^{\dagger}=D^+$ in the space spanned by the genuine eigenfunctions.

It must be emphasized  that $\Psi_z$ is {\it not} an  eigenfunction of
$D$. Indeed, one has

\begin{eqnarray}
D\Psi_z = -D^+ \Psi_z + 2V\Psi_z = z\Psi_z + 2V\Psi_z\per .
\end{eqnarray}

\noindent  This is in accordance with the   analogy
with the  coherent states
which are eigenstates of annihilation operator but not those
of creation operator.

\section{$D^+=D^{\dagger}$ condition and Hermitian form}

The requirement that $D^{+}$ is indeed  the Hermitian conjugate
of $D$ implies that the Hermitian form
satisfies

\begin{eqnarray}
\langle f\vert D^{+} g\rangle   = \langle Df\vert  g \rangle \per . 
\end{eqnarray}

\noindent This condition implies

\begin{eqnarray}
\langle \Psi_{z_1}\vert D^+\Psi_{z_2}\rangle
= \langle D\Psi_{z_1}\vert \Psi_{z_2}\rangle
\per .
\end{eqnarray}

The first (not quite correct) guess is that the Hermitian form 
is defined as an  integral of the product $\overline{\Psi}_{z_1}\Psi_{z_2}$
of the eigenfunctions of the operator $D$ over 
the non-negative real axis using a suitable
integration measure. 
 The Hermitian form  can be defined by continuing
the integrand from the non-negative
real axis to the entire  complex $t$-plane  and noticing that it has a cut
along the non-negative real axis. This suggests  the definition
of the Hermitian form, not as a mere integral over the non-negative real
axis, but as a contour integral along   curve $C$ defined so that
it encloses the non-negative real axis, that is $C$

a) traverses the non-negative real axis along the line $Im[t]=0_-$
from  $t=\infty + i0_-$  to $t=0_+ +i0_-$, 

b) encircles the origin around a small circle
from $t=0_+ +i0_-$   to $t= 0_+ +i0_+$,

c) traverses the  non-negative real axis
along the line $Im[t]=0_+$  from $t=0_+ +i0_+$ to 
$t= \infty+i0_+$ . \\
Here $0_{\pm}$  signifies taking the limit $x= \pm\epsilon$, $\epsilon>0$, $\epsilon \rightarrow 0$.  

 $C$ is the correct choice if   the integrand  defining the
inner product approaches zero sufficiently fast at
the limit  $Re[t]\rightarrow \infty$.  Otherwise
one must assume  that    the integration contour 
 continues along the  circle $S_R$ of radius $R\rightarrow \infty$ 
back to $t=\infty+i0_-$ to form a closed contour. It however turns
out that this is not necessary. 
One can deform the integration contour 
rather  freely: the only constraint is 
 that the deformed integration contour  does not cross over any cut or pole
associated with the analytic continuation of   the integrand
 from the non-negative real axis to the entire
complex plane. 

\vm

Scaling invariance  dictates the form of the  integration
measure appearing in the Hermitian form  uniquely to be $dt/t$.
The Hermitian form  thus obtained also makes  possible
to satisfy the crucial $D^+ = D^{\dagger}$
condition. The Hermitian form is thus defined as

\begin{eqnarray}
\langle f\vert g\rangle   =  -\frac{K}{2\pi i}
\int_C \overline{f}g \frac{dt}{t}\per .
\end{eqnarray}

\noindent      $K$ is a numerical
constant to be determined later. The possibility to
deform the shape of $C$ in wide limits realizes conformal invariance
stating that the change of the shape of the integration contour
induced by a conformal transformation, which is  nonsingular 
inside the integration  contour, leaves the
value of  the contour integral of an analytic function unchanged.
This scaling invariant Hermitian form
is indeed a correct guess. By applying  partial integration one can write

\begin{eqnarray}
\langle \Psi_{z_1}\vert D^+\Psi_{z_2}\rangle
= \langle D\Psi_{z_1}\vert \Psi_{z_2}\rangle - \frac{K}{2\pi i}
\int_C dt \frac{d}{dt} \left[\overline{\Psi}_{z_1}(t) \Psi_{z_2}(t)\right]
\per .
\end{eqnarray}

\noindent  The integral of a total differential comes from the operator 
$L_0=td/dt$ and must vanish.  For a  non-closed  integration contour $C$
the  boundary terms from the partial integration could spoil the  $D^+=D^{\dagger}$ condition unless the eigenfunctions vanish at the end points of the integration contour ($t=\infty + i0_{\pm}$).   

The explicit expression of the Hermitian form is  given by

\begin{eqnarray}
\langle \Psi_{z_1}\vert \Psi_{z_2}\rangle
&=&  -\frac{K}{2\pi i} \int_{C}  \frac{dt}{t} F^2(t) t^{z_{12}}   \per ,\nonumber\\
z_{12}&=&\overline{z}_1+z_2\per .
\label{quadr} 
\end{eqnarray}

\noindent It must be emphasized that it is  
$\overline{\Psi}_{z_1}\Psi_{z_2}$ rather than
eigenfunctions which is continued from the 
non-negative real axis to the  complex $t$-plane: therefore one
indeed obtains an analytic function as a result.
 
 An essential role in the argument claimed to prove the
Riemann hypothesis  is played by 
the crossing symmetry

\begin{eqnarray}
\langle\Psi_{z_1}\vert \Psi_{z_2} \rangle  & =&  
  \langle\Psi_0\vert \Psi_{\overline{z}_1+z_2} \rangle  
\label{cross}
\end{eqnarray}

\noindent of the Hermitian form.  This symmetry is analogous to 
the crossing symmetry of particle physics  stating that 
 the S-matrix is symmetric with respect to the  
replacement  of the particles in the  initial state with
their antiparticles in  the final state or vice versa \cite{field}.

The Hermiticity of the Hermitian form   implies

\begin{eqnarray}
\langle \Psi_{z_1}\vert \Psi_{z_2}\rangle
&=&  
\overline{\langle \Psi_{z_2}\vert \Psi_{z_1}\rangle}\per .
\label{herm} 
\end{eqnarray}

\noindent This condition, which is {\it not} trivially satisfied,
  in fact determines the eigenvalue spectrum.

\section{How to choose the function $F$?}

 The remaining task is to choose the function $F$ in such a manner
that the orthogonality conditions for the solutions $\Psi_0$ and $\Psi_z$
reduce to the condition that $\zeta$  or some
function proportional to $\zeta$ vanishes at the point $-z$. 
The definition of   $\zeta$   based
on   analytical continuation performed 
by Riemann suggests how to proceed.
Recall that the expression of  $\zeta$  converging in
the region $Re[s]>1$ reads \cite{Titchmarch} as

\begin{eqnarray}
\Gamma (s) \zeta (s)= \int_{0}^{\infty} \frac{dt}{t} \frac{exp(-t)}{\left[1-exp(-t)\right]} t^s 
\per . 
\end{eqnarray}
 
\noindent One can analytically continue this expression to
a function defined  in  the entire complex plane by 
noticing that the  integrand is discontinuous along the  cut 
extending from $t=0$ to $t=\infty$. 
Following Riemann it is however more convenient to consider the
discontinuity for a function obtained by multiplying the integrand
with  the factor

 $$(-1)^s \equiv  exp(-i\pi s)  \per .$$

 The discontinuity  $Disc (f)\equiv f(t)-f(texp(i2\pi))$  of the
resulting function is given by

\begin{eqnarray}
Disc\left[\frac{exp(-t)}{\left[1-exp(-t)\right]} (-t)^{s-1}\right]= 
-2i sin(i\pi s)  \frac{exp(-t)}{\left[1-exp(-t)\right]} t^{s-1} \per .
\end{eqnarray}

\noindent  The discontinuity 
 vanishes  at the limit $t\rightarrow 0$ for $Re[s]>1$.
 Hence one can  define  $\zeta$ by modifying the integration
contour from the non-negative  real axis to an integration
contour $C$ enclosing non-negative real axis defined in the
previous section.

This amounts  to writing 
the analytical continuation of $\zeta (s)$ in the form  

\begin{eqnarray}
-2i\Gamma (s) \zeta (s) sin(i\pi s)= 
\int_{C} \frac{dt}{t} \frac{exp(-t)}{\left[1-exp(-t)\right]} 
(-t)^{s-1} \per . 
\label{riemann}
\end{eqnarray}

\noindent This expression equals to $\zeta (s)$  for $Re[s]>1$
and defines $\zeta(s)$ in the entire complex plane since 
the integral around the origin  eliminates the singularity.

The crucial observation is that the 
integrand on the righthand side of Eq.  \ref{riemann}
has precisely the same general form as that appearing
in  the Hermitian form defined  in Eq. \ref{quadr} defined 
using the same integration contour $C$. 
The integration measure is $dt/t$, 
the factor $t^{s}$ is  of the same form as 
the factor $t^{\overline{z}_1 +z_2 }$ appearing 
in the Hermitian form,  and the function $F^2(t)$ is given by

 $$F^2(t)= \frac{exp(-t) }  {1-exp(-t) } \per .$$

\noindent Therefore one can make the identification

\begin{eqnarray}
F(t) =  \left[ \frac{exp(-t)}{1-exp(-t)}\right]^{1/2} \per . 
\end{eqnarray}

\noindent Note that the argument of
the  square root is non-negative on the non-negative
real axis and that $F(t)$
decays exponentially on the non-negative real axis and
has $1/\sqrt{t}$ type singularity at origin. From this it
follows that the eigenfunctions $\Psi_z (t)$ approach zero exponentially
at the limit $Re[t]\rightarrow \infty$ so that one can use the
non-closed integration contour $C$.

With this assumption,  the Hermitian form reduces to the  expression

\begin{eqnarray}
\langle\Psi_{z_1}\vert \Psi_{z_2} \rangle   &=& -\frac{K}{2\pi i }  
\int_{C} \frac{dt}{t} \frac{exp(-t)}{\left[1-exp(-t\right]} 
(-t)^{z_{12}}\nonumber\\
\nonumber\\
&=& \frac{K}{\pi } sin(i\pi z_{12})   \Gamma (z_{12})  \zeta (z_{12}) \per .   
\label{inner}
\end{eqnarray}

\noindent Recall that the definition $z_{12}= \overline{z}_1+z_2$ is
adopted.    Thus the orthogonality of the eigenfunctions is equivalent
to the vanishing of  $\zeta (z_{12})$.

 \section{Study of the Hermiticity condition}

In order to derive information about the
spectrum one must explicitely study what the statement
that $D^{\dagger}$ is Hermitian conjugate of $D$ means.
The defining equation is just the generalization
of the  equation

\begin{eqnarray}
A^{\dagger}_{mn}= \overline{A}_{nm}\per .
\end{eqnarray}

\noindent defining the notion of Hermiticity for matrices. 
Now indices $m$ and $n$ correspond to the eigenfunctions
$\Psi_{z_i}$,  and one obtains

$$\begin{array}{l}
\langle \Psi_{z_1}\vert D^+\Psi_{z_2}\rangle=z_2\langle \Psi_{z_1}\vert D^+\Psi_{z_2}\rangle
= \overline{\langle \Psi_{z_2}\vert D\Psi_{z_1}\rangle }
=\overline{\langle D^+\Psi_{z_2}\vert \Psi_{z_1}\rangle}
= z_2\overline{\langle \Psi_{z_2}\vert \Psi_{z_1}\rangle}\per .\\
\end{array}$$

\noindent   Thus  one has

\begin{eqnarray}
G(z_{12}) &=& \overline{G(z_{21})}= \overline{G(\overline{z}_{12})}\nonumber\\
G(z_{12}) &\equiv& \langle \Psi_{z_1}\vert \Psi_{z_2}\rangle
\per . 
\end{eqnarray}

\noindent  The condition  states that the Hermitian form  defined
by the contour integral is
indeed Hermitian. This is {\it not} trivially true.  
Hermiticity condition obviously determines the spectrum of 
the eigenvalues of $D^+$.

To see the implications of the  Hermiticity condition,
one must study the behaviour of the function $G(z_{12})$ 
under complex conjugation
of both the argument and the value of the function itself. To achieve this
one must write the integral 

$$G(z_{12})   = -\frac{K}{2\pi i } 
\int_{C} \frac{dt}{t} \frac{exp(-t)}{\left[1-exp(-t)\right]} (-t)^{z_{12}}$$

\noindent  in a  form
from which one can easily deduce the behaviour of this function
under complex conjugation. To achieve this, one must
perform the change $t\rightarrow u=log(exp(-i\pi )t)$ of the integration variable
giving

\begin{eqnarray}
G(z_{12})   &=& -\frac{K}{2\pi i}  
\int_{D} du \frac{exp(-exp(u))  }{\left[1-exp(-(exp(u)))\right]} 
exp(z_{12}u) \per . \nonumber\\
\end{eqnarray}

\noindent  Here $D$ denotes the image of the integration contour
$C$ under $t \rightarrow u= log(-t)$.  $D$
is a fork-like  contour which\\
a)  traverses the line  $Im[u]=i\pi$ from $u=\infty+i\pi$ to $u=-\infty +i\pi$ , \\
 b) continues from $-\infty +i\pi$ to $-\infty -i\pi$ along the imaginary
$u$-axis (it is easy to see that the contribution from this part
of the contour vanishes),\\
c) traverses the real $u$-axis from $u=-\infty-i\pi$ to $u=\infty-i\pi$,\\

The integrand differs on the line $Im[u]=\pm i\pi$ from
that on the line $Im[u]=0$ by the factor  $ exp(\mp i\pi z_{12})$ so that
one can write $G(z_{12})$ as integral over real $u$-axis

\begin{eqnarray}
G(z_{12})   &=& -\frac{2K}{\pi }  
sin(i\pi z_{12}) \int_{-\infty}^{\infty} du \frac{exp(-exp(u))  }{\left[1-exp(-(exp(u)))\right]} 
exp(z_{12}u) \per . \nonumber\\
\end{eqnarray}

\noindent 
From this form the effect of the
transformation $G(z)\rightarrow \overline{G(\overline{z})}$
can be deduced. Since the integral is along the real $u$-axis, complex
conjugation amounts only to the replacement  $z_{21}\rightarrow z_{12}$,
and one has

\begin{eqnarray}
  \overline{G(\overline{z}_{12})}    &=& -\frac{2K}{\pi } \overline{sin(i\pi z_{12})} 
\int_{-\infty}^{\infty} du \frac{exp(-exp(u))  }{\left[1-exp(-(exp(u)))\right]} 
exp(z_{12}u) \per  \nonumber\\
&=& -\frac{\overline{sin(\pi z_{12})}}{sin(\pi z_{12})}  G(z_{12})    \per .
\end{eqnarray}

\noindent The substitution of this result  to the Hermiticity condition
gives

\begin{eqnarray}
G(z_{12}) =   -\frac{\overline{sin(i\pi z_{12})}}{sin(i\pi z_{12})} G(z_{12}) \per \mbox{for}\per
x_1+x_2\leq 1\per .
\end{eqnarray}

There are two   manners to satisfy the Hermiticity condition.

a) The condition 
 
\begin{eqnarray}
G(z_{12}) = 0  
\end{eqnarray}

\noindent is the only  manner to satisfy the  Hermiticity
condition for $x_1+x_2<1$ and $y_{2}-y_1\neq 0$.  This implies 
the vanishing of $\zeta$:

\begin{eqnarray}
\zeta (z_{12}) = 0 \per \mbox{for}\per 0<x_1+x_2<1\per ,\per y_1\neq y_2  \per .
\end{eqnarray}

\noindent In particular, this condition
must be  true for $z_1=0$ and $z_2=1/2+iy$.
Hence  the  eigenfunctions with
the  eigenvalue $z=1/2+iy$ correspond to the zeros
of $\zeta$.

b) The condition 

\begin{eqnarray}
\frac{\overline{sin(i\pi z_{12})}}{sin(i\pi z_{12})}=-1\per , 
\end{eqnarray}

\noindent implying

\begin{eqnarray}
 exp(-\pi i (x_1+x_2)) =1\per ,
\end{eqnarray}

\noindent is satisfied. This condition is satisfied for 
$x_1+x_2=n$.  The
highly non-trivial implication is that  the 
states $\Psi_z$ having  {\it real norm}  and  $0<Re[z]<1$
correspond to the zeros of $\zeta$  on the 
line $Re[s]=1/2$. Thus the study of mere Hermiticity
 conditions almost proves the Riemann hypothesis.

\section{Does the  Hermitian form define inner product?}

Before considering the question whether the Hermitian form
defines a positive definite Hilbert space inner product, 
a couple of  comments concerning the general properties
of the Hermitian form  are in order.

a)  The Hermitian form is  proportional to the factor

$$sin(i\pi (y_2-y_1))  \per , $$

\noindent    which vanishes for $y_1=y_2$.
For $y_1=y_2$ and  $x_1+x_2=1$  ($x_1+x_2=0$)   the 
diverging factor $\zeta (1)$ ($\zeta (0)$) compensates
the vanishing of this factor.
Therefore  the  norms of the  eigenfunctions $\Psi_{z}$ with 
 $z=1/2+iy$    must be calculated explicitly from the
defining integral. Since the  contribution from the cut
vanishes in this case, 
one obtains only an integral along a small circle  around 
the origin. This gives   the result

\begin{eqnarray}\begin{array}{ll}
\langle \Psi_{z_1}\vert \Psi_{z_1}\rangle
=K \per \mbox{for} \per z_1=\frac{1}{2}+iy\per , &
\langle \Psi_{0}\vert \Psi_{0}\rangle
=-\frac{K}{2}\per .\\
\end{array}
\end{eqnarray}


\noindent   Thus the norms of the eigenfunctions  are finite.
For $K=1$ the norms of $z=1/2+iy$ eigenfunctions are equal to one.
$\Psi_0$ has however negative norm $-1/2$ so that
the Hermitian form in question is not a genuine  inner product in
the space containing  $\Psi_0$.

b)  For  $x_1=x_2=1/2$ and $y_1\neq y_2$ the factor is nonvanishing and one has

\begin{eqnarray}
\langle \Psi_{z_1}\vert \Psi_{z_2}\rangle
=-\frac{1}{\pi }\zeta (1+i(y_2-y_1))\Gamma (1+i(y_2-y_1)) sinh(\pi (y_2-y_1)) 
\per .\nonumber\\
\end{eqnarray}

\noindent   The nontrivial zeros of $\zeta$ are known to  
belong to the critical strip defined by  $0<Re[s]<1$.  Indeed,
the theorem of Hadamard and de la Vallee  Poussin \cite{analysis}
states the non-vanishing of $\zeta$ on the line $Re[s]=1$.
Since  the  non-trivial zeros of $\zeta$ are located symmetrically
with respect to the line $Re[s]=1/2$,
 this implies that the line $Re[s]=0$ cannot contain
 zeros of $\zeta$.  This result implies  that  the states 
$\Psi_{z=1/2+y}$ are non-orthogonal unless $\Gamma(1+i(y_2-y_1))$
vanishes for some pair of eigenfunctions.

\vm

  It is quite possible that the Hermitian form in question
defines an  inner product in the space spanned by the
  states  $\Psi_z$, $z= 1/2+iy$ having real and positive
norm.  Besides Hermiticity, a necessary condition for this is

$$\vert \langle \Psi_{z_1}\vert \Psi_{z_2}\rangle\vert \leq 1$$

\noindent and gives 

\begin{eqnarray}
-\frac{1}{\pi }\vert \zeta (1+iy_{12})\vert \times\vert \Gamma (1+iy_{12}) 
\times \vert sin(i\pi y_{12})\vert \leq 1 \per ,
\end{eqnarray}

\noindent where the shorthand notation $y_{12}=y_2-y_1$ has been used.
The diagonalized metric is positive definite if 
$G(1/2+iy_{12})$ approaches zero sufficiently fast for large values of
argument $y_2-y_1$ so that the nondiagonal part of the metric
can be regarded as a small perturbation. On physical grounds this is
to be expected since coherent states should have overlap which 
is essentially Gaussian function of   the distance $y_2-y_1$.   
$sin(i\pi y_{12})$ however increases
exponentially and this growth must be compensated by the behaviour
of the  the remaining terms. 

\vm

To get some grasp on the behaviour of the 
 Hermitian metric,  one can use
the integral formula

\begin{eqnarray}
\zeta (s) &=&1 + \frac{1}{s-1} + \frac{1}{\Gamma (s)} I(s)\per , \nonumber\\
I(s)&=& \int_0^{\infty} \frac{dt}{t} \left[\frac{1}{\left[1-exp(-t)\right]} -\frac{1}{t}\right]  exp(-t) t^s \per ,\nonumber\\
\Gamma (s) &= & \int_0^{\infty}  \frac{dt}{t}  exp(-t) t^{s} 
\label{expression}
\end{eqnarray}

\noindent proved already by Riemann. Applying the formula in present case,
one has

\begin{eqnarray}
G(1+iy_{12}) &=& -\frac{1}{\pi i}sinh(\pi y_{12})  \Gamma (1+iy_{12}) 
 (1 + \frac{1}{iy_{12}})\nonumber\\
&-& \frac{1}{\pi i}sinh(\pi y_{12}) I(1+iy_{12})  \per .
\end{eqnarray}

\noindent  Hyperbolic sine increases exponentially as a function of 
$y_{12}$ but cannot  spoil the 
Gaussian decay suggested by the analogy with
the coherent states. One can try ot  demonstrate the  Gaussian behaviour
by  an approximate  evaluation of   the integrals
appearing on the left hand side by changing the integration variable
to  $t=exp(u)$. This gives

\begin{eqnarray}
I(1+iy_{12})&=& \int_{-\infty}^{\infty}  du \left[\frac{1}{\left[1-exp(-e^u))\right]} -\frac{1}{e^u}
\right]  exp(-e^u+u +iy_{12}u) \per ,\nonumber\\
\Gamma (1+iy_{12}) &= & \int_{-\infty}^{\infty}  du  exp(-e^u +u +iy_{12}) \per .
\end{eqnarray}

\noindent The exponential term   has a maximum at $u=0$  and 
vanishes extremely rapidly   as a function of $u$ for $u>0$.
The troublesome feature is that for $u<0$ the integrand decays
only exponentially and  is expected to give 
 a  slowly decreasing contribution  which oscillates as a function
of $y$.  Thus it seems that a  Gaussian, and even an exponential,
 overall decay is excluded.

The analogy with the coherent states  however requires
that the integral   decomposes to a Gaussian term
plus an oscillating remainder   which becomes very small for 
$y=y_{12}$.  Even the inner product property requires that the oscillating
term decays faster than $exp(-\pi y_{12})$ 
as a function of $y_{12}$.
The needed faster than $exp(-\pi y_{12})$ decay
requires that  that the points $y=y_2-y_1$ are approximate
zeros of $G(1+iy)$,  that is approximate zeros of 
 $\zeta (1+iy)$ or, less probably, those  of  $\Gamma (1+iy)$.  
 The mechanism giving rise to an  approximate
zero would be a cancellation of the terms
proportional to  $\Gamma (1+iy)$
and $I(1+iy)$ in the expression of Eq. \ref{expression}
for $\zeta$.  An extremely intricate organization of 
the apparently chaotically located zeros and almost-zeros of $\zeta$ 
is required to guarantee that the Hermitian form defines an inner product.
 Whether the
differences $y=y_2-y_1$ represent approximate
zeros of  $\zeta(1+iy)$ on the  line $Re[s]=1$,  
can be tested numerically.

That the behaviour of $\zeta (1+iy)$  as a function of $y$
can be regarded as a  superposition of a Gaussian term and an  oscillating
term,    is suggested by the following argument. 
If the Gaussian approximation around the origin

$$ -exp (u)+u \simeq -\frac{u^2}{2}\per $$

\noindent were a good approximation,   the integrals in question 
would reduce to Gaussian integrals

\begin{eqnarray}
I(1+iy)&\simeq& \frac{1}{e-1}J(y)\per , \nonumber\\
\Gamma (1+iy) &\simeq &J(y)\per , \nonumber\\
J(y)&=& \int_{-\infty}^{\infty}  du   exp(-\frac{u^2}{2} +iyu) 
= \sqrt{2\pi}exp(-\frac{y^2}{2})
\per .
\end{eqnarray}

\noindent Thus one would have

\begin{eqnarray}
G(1+iy) &\simeq & -\frac{\sqrt{2}}{\sqrt{\pi} i}sinh(\pi y)
exp(-\frac{y^2}{2}) \left[  \frac{1}{iy}+ \frac{e}{e-1}\right]\per .
\end{eqnarray}

\noindent The behaviour would  
be indeed Gaussian for large values of $y$.

Possible problems are also  caused
by the small values of $y_{12}$ for which one might have $\vert G(1+iy{12})
\vert >1$ implying the failure of the Schwartz inequality

\begin{eqnarray}
\vert \langle \Psi_{z_1}\vert \Psi_{z_2}\rangle \vert &\leq& 
\vert \Psi_{z_1}\vert \vert \Psi_{z_2} \vert 
\end{eqnarray}

\noindent chararacterizing positive definite metric.  
In the Gaussian approximation the value of $\vert G(1+iy{12})
\vert$ at the limit $y_{12}=0$ is $\sqrt{2\pi}\simeq 2.5066$ so
that the danger is real.  The direct calculation of 
 $G(1+iy)$ at the limit $y\rightarrow 0$ by using
$\zeta (1+iy)\simeq 1/iy$ however gives 

\begin{eqnarray}
G(1)=1 \per .
\end{eqnarray}

\noindent By a straighforward calculation one can also
verify that $z=1$ is  a local maximum of $\vert G(z)\vert$.

Intuitively it seems obvious that Schwartz inequality must
hold true quite generally.
The point is that the might-be inner product for the superpositions
$\sum_y f(y) \Psi_{1/2+iy}$ and $\sum_y g(y)\Psi_{1/2+iy}$
of $\Psi_z$   describes  net  correlation for the  functions 
$\overline{f}(y)$ and $g(y)$. 
This correlation can be written as

\begin{eqnarray}
\langle f\vert g\rangle&=& \sum_{y_1,y_2} 
\overline{f}(y_1) G(1+i(y_1-y_2))g(y_2)\per .
\end{eqnarray}

\noindent  Since $G(1+i(y_1-y_2))$ decays like Gaussian, 
the correlation of the functions
 $\overline{f}$ and $g$ is determined mainly  
by the correlation   $\overline{f}$ and $g$ 
at  very small distances  $y_1-y_2$. It is obvious that correlation
is largest when $f$ and $g$ resemble each other
maximally, that is when one has $f=g$.

It is easy to see that
arbitrary small values of $y_{12}$ are  unavoidable. 
 The estimate of Riemann  for the number
of the zeros of $\zeta$ in the interval  $Im[s]\in[0,T]$ along 
the line $Re[s]=1/2$ reads as

\begin{eqnarray}
N(T)&\simeq& \frac{T}{2\pi} \left[ log( \frac{T}{2\pi})-1  \right]\per ,
\end{eqnarray}

\noindent and allows to  estimate  the average density $dN_T/dy$ 
of the zeros and to deduce an   upper limit for
the minimum distance $y_{12}^{min}$
between two zeros in the interval  $T$:

\begin{eqnarray}
\frac{dN_T}{dy}&\simeq& \frac{1}{2\pi} \left[ log( \frac{T}{2\pi})-1  \right]\per ,
\nonumber\\
y_{12}^{min}&\leq & \frac{1}{\frac{dN_T}{dy}}= \frac{2\pi }{\left[ log( \frac{T}{2\pi})-1  \right]}\rightarrow 0 \per \mbox{for}\per T\rightarrow \infty 
\per .
\end{eqnarray}

\noindent This implies that arbitrary small values of $y_{12}$ are  unavoidable. 
Thus a rigorous proof for  $\vert G(1+i(y_1+y_2)\vert <1$ for $y_1+y_2\neq 0$
is required.

\section{Superconformal symmetry}

The reduction ad absurdum argument to be discussed below
relies on the assumption that the orthogonality of $\Psi_w$
with $\Psi_0$  for $Re[w]<1/2$ implies the orthogonality
of $\Psi_w$ with all eigenfunctions $\Psi_z$, 
$z=1/2+iy$  zero of $\zeta$. 
In other words, the vanishing of $\zeta (w)$ 
implies the  vanishing of $\zeta (w+z)$ for any zero  $z$ of $\zeta$,
and one has an infinite number of zeros on the line
$Re[s]= Re[w]+1/2$.

This means  the decomposition of
the space   of the eigenfunctions  orthogonal with respect to $\Psi_0$
to a direct sum $V=\oplus_{x<1/2} V_{x} \oplus H_{1/2}$,
such that 
$V_x$  (for which Hermitian form is not inner product)
 contains the  non-orthogonal eigenfunctions $\Psi_{x+iy}$ and $\Psi_{1-x-iy}$ 
and   the spaces $H_x$ and $H_{1/2}$ are orthogonal to each other
for each value of $x$.   The requirement that the eigenfunctions
 having  a positive
norm are orthogonal to the eigenfunctions with  complex norm
and  orthogonal to the state $\Psi_0$,  looks very natural but
it is not easy to justify rigorously  this assumption without
assuming some kind of a symmetry.

Here  superconformal symmetry, which stimulated
the idea behind  the proposed proof of the Riemann hypothesis,
could  come in rescue.
First of all, one can 'understand' the restriction of the non-trivial
zeros to  the line
$Re[s]=1/2$ by noticing that $x$ can be  interpreted  as the real part
of conformal weight defined as eigenvalue of the scaling operator
$L_0= td/dt$ in superconformal field theories \cite{sconf,TGD,octo}. For
the  generators 
of the superconformal algebra,  conformal weights are indeed half-integer
valued. The following construction is essentially a construction
of a second-quantized superconformal quantum field theory for the
system described by $D^+$.  

\vm

One can indeed identify a conformal
algebra naturally associated with the proposed dynamical system. 
 The generators  

\begin{eqnarray}
L_z= \Psi_z D^+
\end{eqnarray}

\noindent generate conformal algebra
with commutation relations ($[A,B]\equiv AB-BA$)

\begin{eqnarray}
\left[L_{z_1},L_{z_2}\right]&=& (z_2-z_1)L_{z_1+z_2}
 \per .
\end{eqnarray}

\noindent The extension of this algebra to superconformal
algebra requires the introduction of the fermionic generators
$G_z$ and $G_z^{\dagger}$. To avoid confusions it must be emphasized
that following convention concerning Hermitian conjugation is adopted
to make notation more fluent:

\begin{eqnarray}
(O_w)^{\dagger}= O^{\dagger}_{\overline{w}}\per .  
\end{eqnarray}

\noindent Fermionic generators $G_z$ and $G_z^{\dagger}$
satisfy the following anticommutation  and commutation
relations:

\begin{eqnarray}\begin{array}{lll}
\{G_{z_1},G^{\dagger}_{z_2}\}= L_{z_1+z_2} \per , & \left[L_{z_1},G_{z_2}\right]=
 z_2G_{z_1+z_2}\per , &  \left[L_{z_1},G^{\dagger}_{z_2}\right]=
 -z_2G^{\dagger}_{z_1+z_2}\per .
\end{array}\nonumber\\
\end{eqnarray}

\noindent This definition differs from that used in the standard approach
\cite{sconf}
in that generators $G_z$ and $G_z^{\dagger}$ are introduced
separately. Usually
one  introduces only the the generators $G_n$ and assumes Hermiticity condition
$G_{-n}= G_n^{\dagger}$.  The anticommutation relations of $G_z$ contain
usually also central extension term. Now this term is not present
as will be found.

Conformal algebras are accompanied by   Kac Moody algebra which
results as a  central extension of  the algebra of 
the local gauge transformations for some Lie group on circle or line
\cite{sconf}.  In the standard approach 
Kac Moody generators are Hermitian in the sense that one has
$T_{-n}=T_n^{\dagger}$ \cite{sconf}. Now this condition is dropped
and one introduces also the generators $T^{\dagger}_z$.
In present case the counterparts for the  generators $T^{\dagger}_z$ of 
the local gauge transformations   act as 
translations $z_1\rightarrow z_1+z$ in   the index
space  labelling eigenfunctions  
and geometrically  correspond to  the  multiplication of $\Psi_{z_1}$
with the  function $t^z$

\begin{eqnarray}
T^{\dagger}_{z_1}\Psi_{z_2}= t^{z_1} \Psi_{z_2} = \Psi_{z_1+z_2} \per .
\end{eqnarray}

\noindent These transformations correspond to the isometries of
the Hermitian form  defined by $G(z_{12})$ and are therefore
natural symmetries at the level of the entire  space
of the eigenfunctions.
 
The commutation relations with the conformal generators
follow from this  definition  and are given by

\begin{eqnarray}\begin{array}{ll}
\left[L_{z_1},T_{z_2}\right]= z_2 T_{z_1+z_2}\per , &
 \left[L_{z_1},T^{\dagger}_{z_2}\right]= 
-z_2 T^{\dagger}_{z_1+z_2}\per ,\end{array}
\end{eqnarray}

\noindent The  central extension making this  
commutative algebra to Kac-Moody algebra
is  proportional to the Hermitian metric

\begin{eqnarray}\begin{array}{lll}
\left[T_{z_1},T_{z_2}\right]=0\per , &\left[T^{\dagger}_{z_1},T^{\dagger}_{z_2}\right]=0\per ,
&\left[T^{\dagger}_{z_1},T_{z_2}\right]=(z_1-z_2)G(z_1+z_2)
 \per .\\\end{array}
\end{eqnarray}

\noindent  One could also consider the central extension $\left[T^{\dagger}_{z_1},T_{z_2}\right]=G(z_1+z_2)$,
which is however  not  the  standard Kac-Moody
central extension.

  One can extend  Kac Moody  algebra to a super Kac Moody algebra
 by adding the fermionic generators $Q_{z}$ and $Q^{\dagger}_z$ 
obeying the anticommutation relations ($\{A,B\}\equiv AB+BA$)

\begin{eqnarray}\begin{array}{lll}
\{Q_{z_1},Q_{z_2}\}= 0\per , & 
\{Q^{\dagger}_{z_1},Q^{\dagger}_{z_2}\}= 0 \per ,
\{Q_{z_1},Q^{\dagger}_{z_2}\}= G(z_1+z_2)\per .
\\
\end{array}
 \end{eqnarray}

\noindent  Note that also  $Q_0$ has a  Hermitian conjugate 
$Q_0^{\dagger}$,   and one has 

\begin{eqnarray}
\{Q_{0},Q^{\dagger}_{0}\}&=& G(0)= -\frac{1}{2}\per 
\end{eqnarray}

\noindent implying that also the fermionic counterpart
of $\Psi_0$ has negative norm.  One can identify  the fermionic
generators as the   gamma matrices of the infinite-dimensional
Hermitian space spanned by the eigenfunctions $\Psi_z$.
By their very definition,  the complexified gamma matrices
$\Gamma_{\bar{z}_1}$ and $\Gamma_{z_2}$
anticommute  to the  Hermitian metric  
$\langle \Psi_{z_1}\vert \Psi_{z_2}\rangle =G(\overline{z}_1+z_2)$.

The commutation relations of the  conformal and  Kac Moody
generators with the   fermionic generators are given by

\begin{eqnarray}\begin{array}{ll}
\left[L_{z_1}, Q_{z_2}\right]= z_2Q_{z_1+z_2}\per , &
\left[L_{z_1}, Q^{\dagger}_{z_2}\right]=  -z_2 Q^{\dagger}_{z_1+z_2}
 \per ,\\
\left[T_{z_1}, Q^{\dagger}_{z_2}\right]=  0\per ,
 &\left[T_{z_1}, Q_{z_2}\right]=  0\per .
\end{array}
\end{eqnarray}

\noindent  The  nonvanishing commutation relations of $T_z$ with $G_z$
and nonvanishing anticomutation
relations of $Q_z$ with 
$G_z$ are given by

\begin{eqnarray}\begin{array}{ll}
\left[G_{z_1}, T^{\dagger}_{z_2}\right]= Q_{z_1+z_2}\per , &
\left[G^{\dagger}_{z_1}, T_{z_2}\right]= -Q^{\dagger}_{z_1+z_2}\per ,\\
\{G_{z_1}, Q^{\dagger}_{z_2}\}= T_{z_1+z_2}\per , &
\{G^{\dagger}_{z_1}, Q_{z_2}\}= T^{\dagger}_{z_1+z_2}\per .
 \end{array}
\end{eqnarray}

\noindent Superconformal generators clearly transform bosonic
and fermionic Super Kac-Moody generators to each other.

 The final step is to construct an explicit representation
for  the generators $G_z$ and $L_z$ in terms of the Super Kac Moody algebra
generators as a generalization of the Sugawara representation
\cite{sconf}. To achieve this, one must introduce the inverse 
 $G^{-1}(z_az_b)$ of the
metric tensor $G(z_az_b)\equiv \langle \Psi_{z_a}\vert \Psi_{z_b}\rangle$,
 which geometrically corresponds to the 
contravariant form of the  Hermitian metric defined by $G$.
Adopting these notations,  one can write the 
generalization for the Sugawara representation
of the superconformal generators as

\begin{eqnarray}
G_z&=& \sum_{z_a} T_{z+z_a} G^{z_az_b}Q^{\dagger}_{z_b} \per , \nonumber\\
G_z^{\dagger}&=& \sum_{z_a} T^{\dagger}_{z+z_a} G^{z_az_b}Q_{z_b} \per . 
\end{eqnarray}

\noindent  One can easily verify that the commutation and anticommutation
relations with the super Kac-Moody generators are indeed correct. 
The generators $L_z$ are obtained as the  anticommutators
of the generators $G_z$ and $G_z^{\dagger}$. Due to
the introduction of the  generators $T_z$, $T^{\dagger}_z$ and
$G_z$, $G^{\dagger}_z$,   the anticommutators $\{G_{z_1}, G^{\dagger}_{z_2}\}$
do not contain any central extension terms. The  expressions for 
the anticommutators  however contains 
terms of form  $T^{\dagger} T Q^{\dagger}Q$ whereas the generators in
the usual Sugawara representation contain only bilinears of
type $T^{\dagger} T$ and $Q^{\dagger} Q$.
The inspiration for introducing  the generators $T_z$,$G_z$
and   $T^{\dagger}_z$,  $G^{\dagger}_z$  separately  comes from the construction
of the physical states as generalized
superconformal representations in  quantum TGD \cite{TGD}.
The proposed algebra differs from the  standard superconformal algebra \cite{sconf}
also in that  the indices $z$ are now complex numbers rather than half-integers or integers as in the case of the ordinary superconformal algebras \cite{sconf}. 
It must be emphasized that one could also consider the   commutation relations $\left[T^{\dagger}_{z_1},T_{z_2}\right]=iG(z_1+z_2)$ and
they might be more the physical choice since $z_2-z_1$ is now a complex
number unlike for ordinary superconformal representations. It is not
however clear how and whether one could construct the counterpart of the Sugawara
representation in this case.

\vm

Imitating the  standard procedure used in
the construction of the representations of the superconformal
algebras \cite{sconf},  one can assume 
that the  vacuum state 
is annihilated   by  {\it all} generators $L_z$
irrespective of the
value of $z$:  

\begin{eqnarray}\begin{array}{ll}
L_z\vert 0\rangle= 0\per , & G_z\vert 0\rangle= 0\per .\\
\end{array}
\end{eqnarray}

\noindent That all generators
 $L_z$ annihilate the vacuum state follows from the representation
$L_z= \Psi_z D_+$ because $D_+$ annihilates $\Psi_0$. If $G_0$
annihilates vacuum then also $G_z\propto [L_z,G_0]$ does the same.

The action of  $T_z^{\dagger}$ on an 
eigenfunction is simply a multiplication by 
$t^z$: therefore one cannot require that $T_z$ annihilates 
the vacuum state as is usually done \cite{sconf}.
The action of $T_0$  is multiplication by $t^0=1$ so that $T^0$
and $T^{\dagger}_0$ act as unit operators in the  space
of the physical states.
In particular, 

\begin{eqnarray}
T_0\vert 0\rangle =T_0^{\dagger}\vert 0\rangle= \vert 0\rangle \per .
\end{eqnarray}

\noindent   This implies the condition 

\begin{eqnarray}
\left[T_0, T_z^{\dagger}\right]= iz G(z)=0
\end{eqnarray}

\noindent in the space of the physical states so that physical states
must correspond to the zeros of $\zeta$
and possibly to   $z=0$.  Thus one can generate the physical states
from vacuum by acting using operators $Q_z^{\dagger}$ and $T^{\dagger}_z$
with $\zeta (z)=0$. If one requires that the  physical states also have
 real and positive
norm squared, only the zeros of $\zeta$ on the  line $Re[s]=1/2$ are allowed.
Hence  the requirement that a unitary representation of
the  superconformal
algebra is in question,  forces Riemann hypothesis.

It is important to notice that 
$T^{\dagger}_z$ and $Q^{\dagger}_z$ cannot annihilate the  vacuum:
this would lead to the condition $G(z_1+z_2)=0$ implying 
the vanishing of $\zeta (z_1+z_2)$ for any pair $z_1+z_2$.
One can however assume that $Q_z$ annhilates the vacuum state

\begin{eqnarray}
Q_z\vert 0\rangle =0 \per .
\end{eqnarray}

This  inspires the hypothesis that only 
the generators with conformal weights  $z=1/2+iy$ 
generate physical states from vacuum   realizable in the space of the
eigenfunctions $\Psi_z$ and their fermionic counterparts.  
This means that the action of the bosonic generators 
$T^{\dagger}_{1/2+iy}$  and  fermionic   generators $Q_0^{\dagger}$
and  $Q^{\dagger}_{1/2+iy}$,  as well as 
the action of the corresponding superconformal generators
$G^{\dagger}_{1/2+iy}$, generates  bosonic and fermionic states with conformal
weight $z=1/2+iy$  from the vacuum state:

\begin{eqnarray}\begin{array}{ll}
\vert 1/2+iy\rangle_B \equiv T^{\dagger}_{1/2+iy}\vert 0\rangle\per ,&
\vert 1/2+iy\rangle_F \equiv Q^{\dagger}_{1/2+iy}\vert 0\rangle\per .
\end{array}
\end{eqnarray}

\noindent One can identify  the states generated by the 
Kac Moody generators $T^{\dagger}_z$ from 
the vacuum as the eigenfunctions $\Psi_z$. 
The system as a whole represents a second quantized supersymmetric version
of the bosonic system defined by the eigenvalue equation
for $D^+$ obtained by assigning to each eigenfunction
a fermionic counterpart and performing second quantization
as a free quantum field theory.

\section{Is the proof of the Riemann hypothesis by reductio ad absurdum
possible using superconformal invariance?}

 Riemann hypothesis is proven if all eigenfunctions for which 
the Riemann
Zeta function vanishes,  correspond to the states having a  real 
 and positive norm squared.   The expectation is that superconformal
invariance realized in some sense excludes all zeros
of $\zeta$ except those on the line $Re[s]=1/2$. 
The problem is to define precisely what one means
with superconformal invariance and one can generate 
large number of reduction
ad absurdum type proofs depending on how superconformal
invariance is assumed to be realized.

\vm

 The most conservative
option is that superconformal invariance is realized in 
the standard sense.   The action
of the ordinary superconformal generators $L_n$, and $G_n$,  $n\neq 0$ 
on the  vacuum states   $\vert 0 \rangle_{B/F} $ or on any state
$\vert 1/2+iy\rangle_{B/F} $  indeed creates zero norm states as is obvious
from the vanishing of the factor $sin(i\pi z_{12})=sin(\pi (x_1+x_2))$
associated with the inner inner products of these states.
Thus the zeros of $\zeta$ define an infinite family of ground states
for the representations of the ordinary superconformal algebra.
A  generalization of this hypothesis is that the action
of $L_n$ and $G_n$, $n\neq 0$,  on any state 
$\vert w\rangle_{B/F} $, $\zeta (w)=0$,  creates   states
which are mutually orthogonal zero norm states. 
This implies $\zeta (n+2Re[w])=0$ for all values of $n\neq 0$ and,
 since the  real axis contains zeros of $\zeta$ only at
the points $Re[s]=-2n$, $n>0$,
leads to a reductio ad absurdum unless one has $Re[w]=1/2$. 
Thus the proof of the Riemann hypothesis would
reduce to showing that the  action of the ordinary 
superconformal algebra generates mutually orthogonal zero norm 
states from any state $\vert w\rangle_{B/F}$ with $\zeta (w)=0$.
The proof of this physically plausible hypothesis is not obvious.

\vm

One can imagine also other strategies. The minimal requirement
is certainly that some subalgebra of the superconformal algebra generates
a space of states satisfying the Hermiticity condition.
   The quantity

\begin{eqnarray}
\Delta (\overline{w}_1 +w_2)&\equiv& 
\langle  w_1\vert w_2\rangle -\overline{\langle  w_2\vert w_1\rangle} 
=  G(\overline{w}_1+w_2)- \overline{G(\overline{w}_2+w_1)}
\end{eqnarray}

\noindent must define the conformal invariant in question  since this
quantity must vanish in the space of the physical states for which
the  metric is Hermitian.  This requirement does not
however imply anything nontrivial for the ordinary conformal
algebra having generators $L_n$ and $G_n$: for $Re[w]\neq 1/2$
the condition is indeed satisfied because $G(n+2Re[w])$ does {\it not}
satisfy the Hermiticity condition for any value of $n$.

One can   try to abstract some property of the states associated
with the zeros of $\zeta$ on the line $Re[s]=1/2$. 
The generators $L_{1/2-iy}$ and  $G_{1/2-iy}$ generate 
zero norm states   from the states $\vert 1/2+iy \rangle_{B/F}$,
when $1/2+iy$ corresponds to the zero of $\zeta$ on the line 
$Re[s]=1/2$.  One can try to generalize this observation
so that it applies to an arbitrary
state  $\vert w\rangle_{B/F}$, $\zeta (w)=0$.
The generators $L_{1-\overline{w}}$ and  $G_{1-\overline{w}}$ 
certainly generate zero norm states from the states $\vert w\rangle_{B/F}$.
Also the Hermiticity condition  holds   true identically 
and does not have nontrivial implications. One can however consider
alternative generalizations by assuming
that 

a)  either   the generators $L_{\overline{w}}$ and  $G_{\overline{w}}$
or 

b)   $L_{1/2+iy}$ and  $G_{1/2+iy}$ 
generate from the states $\vert w \rangle_{B/F}$, $\zeta (w)=0$
 states satisfying the
Hermiticity condition.

These two hypothesis lead to two
 versions of a  reductio ad absurdum
argument. Suppose  that  $w$ is a zero of $\zeta$. 
This means that the inner product of the states $Q_0^{\dagger}\vert 0\rangle$
and $Q_w^{\dagger}\vert 0\rangle$  and thus also   $\Delta (w)$
vanishes:

\begin{eqnarray}
\begin{array}{ll}
  \langle 0\vert  Q_0 Q_w^{\dagger}\vert 0 \rangle=0\per , 
& \Delta (w)=0\per .
\end{array}
\label{invariant1}
\end{eqnarray}

a)  By acting on this  matrix element by the conformal
algebra generator
$L_{\overline{w}}$ (which acts like derivative operator on 
the arguments
of the should-be Hermitian form), 
and using the
fact that $L_{\overline{w}}$ annilates the  vacuum state,
one obtains

\begin{eqnarray}
 \langle 0\vert  Q_0 Q_{\overline{w}+w}^{\dagger}\vert 0\rangle = G(w+\overline{w})
\per .
\label{invariant2}
\end{eqnarray}

\noindent The requirement  $\Delta (w+\overline {w})=0$ 
implies the reality of $G(w+\overline{w})$  and thus   
the condition $Re[w]=1/2$ leading to  the Riemann hypothesis.
Note that  the  argument implying the reality
of $G(w+\overline{w})$  assumes 
only that $L_w$ annihilates vacuum.

If this line of approach is correct,
the basic challenge  would be to show on the  basis of
the  superconformal
invariance alone that the condition  $\zeta (w)=0$
implies that the generators   $L_{\overline{w}}$ 
and $G_{\overline{w}}$ generate new ground states
satisfying the Hermiticity condition.

b) An alternative line of  argument  uses  only the invariance under the
generators $L_{1/2+iy}$ associated with the zeros of $\zeta$, 
and thus certainly belonging to the conformal algebra associated with
the physical states. By 
 applying  the generators $L_{1/2+iy_i}$ to the 
the matrix element $ \langle 0\vert  Q_0 Q_w^{\dagger}\vert 0\rangle=0$
and requiring that Hermiticity is respected,  one can 
 deduce that $G(w+1/2+iy_i)$ satisfies the Hermiticity condition.
Hence the line $Re[s]=Re[w]+1/2$,
and by the  reflection symmetry
also the line $Re[s]=1/2-Re[w]$, contain an infinite number
of zeros of $\zeta$ if one has  $Re[w]\neq 1/2$.
By repeating this process once for the zeros on the line 
$Re[s]=1/2-Re[w]$,  one finds that the lines
$Re[s]=1-Re[w]$ and $Re[s]=Re[w]$   contain infinite number of the  zeros of
$\zeta$  of form $w_{ij}= w+i(y_i+y_j)$, where
$y_i$ and $y_j$ are associated with the  zeros of $\zeta$
on the line $Re[s]=1/2$. By applying this two-step procedure
repeatedly,  one can fill the
lines $Re [s]=  Re[w]$,$1-Re[w]$, $1/2-Re[w]$, $1/2+Re[w]$ 
with the  zeros of $\zeta$.

To sum up, contrary to the original over-optimistic beliefs inspired
by the beauty of the proposed quantum  model,  Riemann hypothesis
demonstrates once again that it is 
equally  resistible against  proof as
it is capable of stimulating new  mathematical ideas.  One
might however hope that superconformal invariance could
in one of the proposed forms or in some other form be
used to rigorously prove Riemann hypothesis.

\vl

\end{document}